\documentclass[a4paper,leqno,12pt]{amsart}
\usepackage{amssymb}
\usepackage{mathrsfs}

\theoremstyle{plain}
\newtheorem{lemma}{Lemma}[section]
\newtheorem{Proposition}[lemma]{Proposition}
\newtheorem{theorem}[lemma]{Theorem}
\newtheorem{coro}[lemma]{Corollary}
\newcommand{\Prop}{\begin{Proposition}}
\newcommand{\enprop}{\end{Proposition}}
\newcommand{\Lemma}{\begin{lemma}}
\newcommand{\enlemma}{\end{lemma}}
\newcommand{\Theorem}{\begin{theorem}}
\newcommand{\entheorem}{\end{theorem}}
\newcommand{\Cor}{\begin{coro}}
\newcommand{\encor}{\end{coro}}

\theoremstyle{remark}
\newtheorem*{remark}{Remark}

\newcommand{\D}{\ensuremath{\mathcal{D}}}
\newcommand{\Z}{\mathbb{Z}}

\renewcommand{\H}{\operatorname{H}}

\renewcommand{\O}{\mathcal{O}}
\newcommand{\SH}{\mathcal{H}}
\newcommand{\g}{\mathfrak{g}}

\newcommand{\Gr}{\operatorname{Gr}}

\newcommand{\To}{\xrightarrow{\phantom{aa}}}
\newcommand{\eq}{\begin{eqnarray}}
\newcommand{\eneq}{\end{eqnarray}}
\newcommand{\eqn}{\begin{eqnarray*}}
\newcommand{\eneqn}{\end{eqnarray*}}
\renewcommand{\phi}{\varphi}

\newcommand{\C}{\mathbb{C}}

\newcommand{\Tor}{\operatorname{Tor}}
\newcommand{\codim}{\operatorname{codim}}

\begin{document}

\title{Local cohomology and \D-affinity\\ in positive characteristic}
\author{Masaki Kashiwara}% \and Niels Lauritzen}
\address{Research Institute for the Mathematical Sciences\\Kyoto University\\Kyoto, Japan}
\email{masaki@kurims.kyoto-u.ac.jp}

\thanks{The research of the first author is partially supported by
Grant-in-Aid for Scientific Research (B) 13440006,
Japan Society for the Promotion of Science. The research of the second
author was supported in part by the TMR-programme ``Algebraic Lie Representations''
(ECM Network Contract No. ERB FMRX-CT 97/0100) and the Danish Natural Science
Research Council.}

\author{Niels Lauritzen}
\address{Institut for matematiske fag\\Aarhus Universitet\\ \AA rhus, Denmark}
\email{niels@imf.au.dk}

\maketitle
\section{Introduction}
\label{sect:intro}
Let $k$ be a field.
Consider
the polynomial ring 
$$
R = k\left[
\begin{array}{ccc}
X_{11} & X_{12} &  X_{13}\\
X_{21} & X_{22} & X_{23}
\end{array}
\right]
$$ 
and $I\subset R$ the ideal generated by the three $2\times 2$ minors
$$
f_1{:=}\left|
\begin{array}{cc} 
X_{12} & X_{13} \\
X_{22} & X_{23}
\end{array}
\right|,\quad
f_2{:=}\left|
\begin{array}{cc}
X_{13} & X_{11}\\
X_{23} & X_{21}
\end{array}
\right|\quad\mbox{and}\quad
f_3{:=}\left|
\begin{array}{cc}
X_{11} & X_{12}\\
X_{21} & X_{22}
\end{array}
\right|.
$$
If $k$ is a field of positive
characteristic, the local cohomology modules $\H^j_I(R)$ 
vanish for $j>2$ (see Chapitre III, Proposition (4.1) in \cite{PS}). 
However, if $k$ is a field of characteristic zero, 
$\H^3_I(R)$ is non-vanishing (see Proposition \ref{prop:nv} 
of this paper or Remark 3.13 in \cite{HL}). 

Consider the Grassmann variety $X = \Gr(2, V)$
of $2$-dimensional vector subspaces of
a $5$-dimensional vector space $V$ over $k$.
Let us take a two dimensional subspace $W$ of $V$.
Then the singularity $R/I$ appears in the Schubert variety 
$Y\subset X$ of $2$-dimensional
subspaces $E$ such that $\dim(E\cap W)\ge1$.
Therefore $\SH^3_Y(\O_X)$ does not vanish
in characteristic zero while it does vanish in positive characteristic.

In this paper we
show how this difference in the vanishing of local cohomology 
translates into 
a non-vanishing first cohomology group for
the $\D_X$-module $\SH^2_Y(\O_X)$ in positive characteristic. 

Previous work of Haastert (\cite{Haastert}) showed 
the Beilinson-Bernstein equivalence (\cite{BB}) to hold 
for projective spaces and the flag manifold 
of $SL_3$
%$\SL_3/B$ 
in positive characteristic.
However, as we show in this paper,
\D-affinity breaks down for 
the flag manifold of $SL_5$
%$X=\SL_5/B$ 
in all positive characteristics.
The Beilinson-Bernstein equivalence, therefore, does not 
carry over to flag manifolds in positive characteristic.

\section{Local cohomology}

Keep the notation from \S \ref{sect:intro}.
A topological proof of the following proposition 
is given in \S \ref{sec:proof}.
\Prop\label{prop:nv}
$H^3_I(R)$ does not vanish in characteristic zero.
\enprop

\Cor $\SH^3_Y(\O_X)$ does not vanish in characteristic zero.
\encor

The local to global spectral sequence
$$
\H^p(X, \SH^q_Y(\O_X))\Rightarrow \H^{p+q}_Y(X, \O_X)
$$
and $\D$-affinity in characteristic zero (\cite{BB}) implies 
$$\H^3_Y(X, \O_X) = \Gamma(X, \SH^3_Y(\O_X))\neq 0.$$ 
On the other hand, if $k$ is a field of positive 
characteristic, $\SH^q_Y(\O_X) = 0$ if $q \neq 2$, since $Y$ is 
a codimension two Cohen Macaulay subvariety of the smooth variety $X$
(\cite{PS}, Chapitre III, Proposition (4.1)). 
This gives a 
totally
%fundamentally 
different degeneration of
the local to global spectral sequence. In the positive characteristic case
we get
$$
\H^p(X, \SH^2_Y(\O_X)) \cong \H^{p+2}_Y(X, \O_X).
$$
We will prove that $\H^3_Y(X, \O_X)\neq 0$ 
even if $k$ is a field of positive characteristic. This will
give the desired non-vanishing
$$
\H^1(X, \SH^2_Y(\O_X)) \neq0
$$
in positive characteristic.

\section{Lifting to $\Z$}

To deduce the non-vanishing of $\H^3_Y(X, \O_X)$ in positive characteristic,
we need to compute the local cohomology over $\Z$ and proceed by
base change. Flag manifolds and their Schubert
varieties admit flat lifts to $\Z$-schemes. In this section $X_\Z$ and
$Y_\Z$ will denote flat lifts of a flag manifold $X$ and a 
Schubert variety $Y\subset X$ respectively.

The local Grothendieck-Cousin complex (cf.~\cite{Kempf}, \S8) of the 
structure sheaf $\O_{X_\Z}$
\begin{equation}\label{eq:res}
 \SH^0_{X_0/X_1}(\O_{X_\Z})\rightarrow \SH^1_{X_1/X_2}(\O_{X_\Z})\rightarrow \cdots, %\tag{$*$}
\end{equation}
where $X_i$ denotes the union of Schubert schemes of codimension $i$, is
a resolution of $\O_{X_\Z}$, since $\O_{X_\Z}$ is Cohen Macaulay, $\codim\, X_i \geq i$
and $X_i\setminus X_{i+1}\rightarrow X$ are affine morphisms for 
all $i$ (see \cite{Kempf}, Theorem 10.9). 
The sheaves in this resolution decompose into direct sums 
$$
\SH_{X_i/X_{i+1}}^i(\O_{X_\Z}) = \bigoplus_{\codim(C) = i} \SH^i_C(\O_{X_\Z})
$$
of  local cohomology  sheaves $\SH^i_C(\O_{X_\Z})$ with support
in Bruhat cells $C$ of codimension $i$. 
The degeneration of the local to global spectral sequence gives
$$
\H^p_{Y_\Z}(X_\Z, \SH^c_C(\O_{X_\Z}))=\H^{p+c}_{C\cap Y_\Z}(X_\Z, \O_{X_\Z}),
$$
since $\SH^i_C(\O_{X_\Z})= 0$ if $i\neq c=\codim(C)$. 

%In fact
Since the scheme $X_i\setminus X_{i+1}$ is affine it follows that
$\H^p_C(X_\Z, \O_{X_\Z})=0$ if $p\neq \codim(C)$ (\cite{Kempf}, Theorem 10.9). 
This shows that the resolution \eqref{eq:res} is acyclic for 
the functor  $\Gamma_{Y_\Z}$ and 
$$
\Gamma_{Y_\Z}(X_\Z, \SH^c_C(\O_{X_\Z})) = 
\left\{
\begin{array}{ll}
0 & \mbox{if $C \not\subset Y_Z$}\\[2pt]
\H^c_C(X_\Z, \O_{X_\Z}) & \mbox{if $C \subset Y_\Z$,}
\end{array}
\right.
$$
where $c = \codim(C)$.
Applying $\Gamma_{Y_\Z}(X_\Z, -)$ to \eqref{eq:res} we get the complex
$$
M^\bullet_Y : \H^c_{C_Y}(X_\Z, \O_{X_\Z}) \rightarrow 
\bigoplus_{\codim(C) = c+1,\,C\subset Y_\Z} \H^{c+1}_{C}(X_\Z, \O_{X_\Z})\rightarrow\cdots,
$$
where $c$ is the codimension of
$Y_\Z$, $C_Y$ is the open Bruhat cell in $Y_\Z$ 
and $\H^c_{C_Y}(X_\Z, \O_{X_\Z})$ sits
in degree $c$. 
Notice that $\H^i(M^\bullet_Y) = 
\H^i_{Y_\Z}(X_\Z, \O_{X_\Z})$ and that $M^\bullet_Y$ is a complex of free abelian
groups. In fact the individual entries $\H^i_C(X_\Z, \O_{X_\Z})$ are
direct sums of weight spaces, 
which are finitely generated free
abelian groups (cf.~\cite{Kempf}, Theorem 13.4). By weight spaces we mean
eigenspaces for a fixed $\Z$-split torus $T$.
The differentials in $M^\bullet_Y$ being $T$-equivariant,
the complex $M^\bullet_Y$ is a direct sum of
complexes of finitely generated free abelian groups.
Since $\H^i(M^\bullet_Y) = \H^i_{Y_\Z}(X_\Z, \O_{X_\Z})$,
one obtains the following lemma.
\Lemma
\label{Lemma:directsum}
Every local cohomology group
$\H^i_{Y_\Z}(X_\Z, \O_{X_\Z})$ is a direct sum of finitely generated
abelian groups. In the codimension $c$ of $Y_\Z$ in $X_\Z$, 
$\H_{Y_\Z}^c(X_\Z, \O_{X_\Z})$ is a free abelian group.
\enlemma

\section{The counterexample}

For a field $k$, let us set
$X_k=X_\Z\otimes k$ and $Y_k=Y_\Z\otimes k$.
Then one has
$H^q_{Y_k}(X_k,\O_{X_k})=\H^q_{Y_\Z}(X_\Z,\O_{X_\Z}\otimes k)$.
Since $\O_{X_\Z}$ is flat over $\Z$, one has
a spectral sequence
\[%E_2^{pq}=
\Tor_{-p}^\Z(\H^q_{Y_\Z}(X_\Z,\O_{X_\Z}),k)\Rightarrow
\H^{p+q}_{Y_k}(X_k,\O_{X_k})\,.
\]
This shows that the natural homomorphism
$$
\H^i_{Y_\Z}(X_\Z,\O_{X_\Z})\otimes k \To \H^i_{Y_k}(X_k,\O_{X_k})
$$
is an injection, 
and it is an isomorphism if the field $k$ is flat over $\Z$. 

In our example (cf.~\S \ref{sect:intro}), one has
$$
\H^3_{Y_\Z}(X_\Z,\O_{X_\Z})\otimes \mathbb{C} 
\cong \H^3_{Y_\C}(X_\C,\O_{X_\C})\neq 0.
$$

By Lemma \ref{Lemma:directsum}, 
the cohomology
%this means that 
$\H^3_{Y_\Z}(X_\Z,\O_{X_\Z})$ 
must contain $\Z$ as a direct summand. Therefore the injection
$$
\H^3_{Y_\Z}(X_\Z,\O_{X_\Z})\otimes k \To
\H^3_{Y_k}(X_k, \O_{X_k}) 
$$
shows that $\H^3_{Y_k}(X_k, \O_{X_k})$ is non-vanishing for any field $k$
of positive characteristic.
Since $\H^3_{Y_k}(X_k, \O_{X_k})
\cong \H^1(X_k, \SH^2_{Y_k}(\O_{X_k}))$, one obtains the following result.

\Prop
$\H^1(X_k, \SH^2_{Y_k}(\O_{X_k}))\neq0$ if $k$ is of positive characteristic.
\enprop

\section{Proof of non-vanishing of $H^3_I(R)$}\label{sec:proof}

In this section, we shall give a topological proof of
Proposition \ref{prop:nv}.
We may assume that the base field is the complex number field $\C$.

The local cohomologies $H^*_I(R)$ are the cohomology groups of the complex
\eqn
&&R\To R[f_1^{-1}]\oplus R[f_2^{-1}]\oplus R[f_3^{-1}]\\
&&\hspace{20pt}
\To R[(f_1f_2)^{-1}]\oplus R[(f_2f_3)^{-1}]\oplus R[(f_1f_3)^{-1}]
\To R[(f_1f_2f_3)^{-1}]\,.
\eneqn
Hence one has

$$H^3_I(R)=\dfrac{R[(f_1f_2f_3)^{-1}]}{R[(f_1f_2)^{-1}]+R[(f_2f_3)^{-1}]
+R[(f_1f_3)^{-1}]}.$$
In order to prove the non-vanishing of $H^3_I(R)$,
it is enough to show
\eq
\dfrac{1}{f_1f_2f_3}\notin 
R[(f_1f_2)^{-1}]+R[(f_2f_3)^{-1}]+R[(f_1f_3)^{-1}]\,.
\label{eq:5.0}
\eneq

Consider the $6$-cycle
\eqn
\gamma&=&\left\{\left(
\begin{array}{ccc}
-t_2u +t_3\bar v& u &  -t_1\bar v\\
-t_2 v-t_3\bar u &v & t_1\bar u
\end{array}\right)\right.
\,\\
&&\left.\hspace{50pt};\,|t_1|=|t_2|=|t_3|=1,\ |u|^2+|v|^2=1
\rule[-1em]{0em}{2em}\right\}\\
&=&\left\{k\left(
\begin{array}{ccc}
-t_2 & 1& 0\\
-t_3& 0& t_1
\end{array}
\right)\right.\\
&&\left.\hspace{50pt};\,|t_1|=|t_2|=|t_3|=1,\ 
k=\left(
\begin{array}{cc}
u &-\bar v\\
v &\bar u
\end{array}\right)
\in SU(2)\rule[-1em]{0em}{2em}\right\}
\eneqn
in $X\setminus (f_1f_2f_3)^{-1}(0)$,
where $X={\mathrm{Spec}}(R)\cong\C^6$.
Then on $\gamma$ one has
$$f_1=t_1,\ f_2=t_1t_2\mbox{ and }f_3=t_3\,.$$
Set $\omega=\bigwedge dX_{ij}$.
Then one has
$\omega=t_1dt_1dt_2dt_3\theta$ on $\gamma$, where
$\theta$ is a non-zero invariant form on $SU(2)$.
%$\{(u,v);|u|^2+|v|^2=1\}$
Therefore one has
$$
\int_\gamma\dfrac{\omega}{f_1f_2f_3}=\int_\gamma\dfrac{dt_1dt_2dt_3\theta}
{t_1t_2t_3}
%=\int_{|t_0|=1}\dfrac{dt_0}{t_0}\int_{|t_0|=1}\dfrac{dt_0}{t_0}
%\int_{|t_0|=1}\dfrac{dt_0}{t_0}\int_{SU(2)}\theta
\not=0.$$
Hence, in order to show \eqref{eq:5.0},
it is enough to prove that %On the other hand, one has
\eq&&
\displaystyle\int_\gamma \phi{\omega}=0
\label{eq:5.1}
\eneq
for any 
$\phi\in R[(f_1f_2)^{-1}]+R[(f_2f_3)^{-1}]+R[(f_1f_3)^{-1}]\,$.

For $\phi\in R[(f_1f_2)^{-1}]$, the equation
%(resp. $\phi\in R[(f_1f_3)^{-1}]$)
\eqref{eq:5.1} holds because we can shrink the cycle $\gamma$
by $|t_3|=\lambda$ from $\lambda=1$ to $\lambda=0$.
For $\phi\in R[(f_1f_3)^{-1}]$),
the equation \eqref{eq:5.1} holds because we can shrink the cycle $\gamma$
by $|t_2|=\lambda$ from $\lambda=1$ to $\lambda=0$.

Let us show \eqref{eq:5.1} for $\phi\in R[(f_2f_3)^{-1}]$.
Let us deform the cycle $\gamma$ by
\eqn
\gamma_{\lambda}&=&\left\{k\left(
\begin{array}{ccc}
-(1-\lambda)t_2&1 & - \lambda t_1t_2t_3^{-1}\\
-t_3& 0& t_1
\end{array}
\right)\,
\right.\\&&\hspace{50pt}\left.
\,;\,|t_1|=|t_2|=|t_3|=1,\ k\in SU(2)\rule[-1em]{0em}{2em}\right\}\,.
\eneqn
Note that the values of $f_1$, $f_2$ and $f_3$
do not change under this deformation.
Hence $\gamma_\lambda$ is a cycle in $X\setminus (f_1f_2f_3)^{-1}(0)$.
One has
%Since 
%$\left(
%\begin{array}{cc}
%t_3&0\\
%0& t_3^{-1}
%\end{array}\right)$
%is in $SU(2)$, one has
\eqn
\gamma_{1}&=&\left\{k
%\left(
%\begin{array}{cc}
%t_3&0\\
%0& t_3^{-1}
%\end{array}\right)
\left(
\begin{array}{ccc}
0&1 & - t_1t_2t_3^{-1}\\
-t_3& 0& t_1
\end{array}
\right)\,;\,|t_1|=|t_2|=|t_3|=1,\ k\in SU(2)\right\}\\
&=&\left\{k
\left(
\begin{array}{ccc}
0&1 &  -t_2\\
-t_3& 0& t_1
\end{array}
\right)\,;\,|t_1|=|t_2|=|t_3|=1,\ k\in SU(2)\right\}\,.
\eneqn
In the last coordinates of $\gamma_1$, one has
$f_2=t_2t_3$ and $f_3=t_3$.
Hence, for $\phi\in R[(f_2f_3)^{-1}]$,
\[\int_{\gamma}\phi\omega=\int_{\gamma_1}\phi\omega\]
vanishes because we can shrink the cycle $\gamma_1$
by $|t_1|=\lambda$ from $\lambda=1$ to $\lambda=0$.

\begin{remark}
Although we do not give a proof here,
$H^3_I(R)$ is isomorphic to $H^6_{J}(R)$ as a D-module.
Here $J$ is the defining ideal of the origin.
\end{remark}

\end{document}